\documentclass[11pt]{article}
\usepackage{fullpage}
\usepackage{latexsym}
\usepackage{epsfig}

\def\qed{$\hfill{\Box}$}

\newfont{\bbold}{msbm10 scaled \magstep1}
\newfont{\bbolds}{msbm7 scaled \magstep1}

\newcommand{\ns}{\mbox{\bbold N}}

\newcommand{\rs}{\mbox{\bbold R}}

\newcommand{\bm}[1]{\mbox{\boldmath \ensuremath{#1}}}

\newcommand{\beq}{\begin{equation}}
\newcommand{\eeq}{\end{equation}}

\newtheorem{Theorem}{Theorem}
\newtheorem{Proposition}{Proposition}
\newtheorem{Corollary}{Corollary}

\title{On a Conjecture of Ira Gessel}
\author{Marko Petkov\v sek \\
Herbert S.\ Wilf}
\date{May 30, 2008}

\begin{document}

\maketitle

\begin{abstract}
Let $F(m;\, n_1,n_2)$ denote the number of lattice walks from $(0,0)$ to $(n_1,n_2)$, always staying in the 
first quadrant $\{(n_1,n_2);\ n_1 \ge 0, n_2 \ge 0\}$ and having exactly $m$ steps, each of which belongs to the set 
$\{E=(1,0),W=(-1,0),NE=(1,1),SW=(-1,-1)\}$. Ira Gessel conjectured that $F(2n;\,0,0) = 
16^n \frac{(1/2)_n(5/6)_n}{(2)_n(5/3)_n}$. We pose similar conjectures for some other values of $(n_1,n_2)$,
and give closed-form formulas for $F(n_1;\,n_1,n_2)$ when $n_1 \ge n_2$ as well as for $F(2n_2 - n_1;\,n_1,n_2)$
when $n_1 \le n_2$. In the main part of the paper, we derive a functional equation satisfied by the generating function
of $F(m;\, n_1,n_2)$, use the kernel method to turn it into an infinite lower-triangular system of linear equations 
satisfied by the values of $F(m;\,n_1,0)$ and $F(m;\,0,n_2) + F(m;\,0,n_2 - 1)$, and
express these values explicitly as determinants of lower-Hessenberg matrices with unit superdiagonals whose non-zero
entries are products of two binomial coefficients.
\end{abstract}

\section{Introduction}

Let $F(m;\, n_1,n_2)$ denote the number of lattice walks from $(0,0)$ to $(n_1,n_2)$, always staying in the 
first quadrant $\{(n_1,n_2);\ n_1 \ge 0, n_2 \ge 0\}$ and having exactly $m$ steps, each of which belongs to the set 
$\{E=(1,0),W=(-1,0),NE=(1,1),SW=(-1,-1)\}$. From the obvious recurrence
\begin{eqnarray}
F(m;\,n_1,n_2) &=& F(m-1;\,n_1+1,n_2)\ +\ F(m-1;\,n_1-1,n_2) \nonumber \\
&+&F(m-1;\,n_1+1,n_2+1)\ +\ F(m-1;\,n_1-1,n_2-1) \label{rec}
\end{eqnarray}
valid for $m\ge 1, n_1,n_2\ge 0$, the initial conditions
\beq
\label{init}
F(0;\,n_1,n_2)\ =\ \left\{\begin{array}{ll}
1, & n_1 = n_2 = 0,\\
0, & {\rm otherwise}
\end{array}\right.
\eeq
(so called because we'll think of $m$ as being the time variable), and the boundary conditions
\beq
\label{bdry}
F(m;\,n_1,n_2)\ =\ 0, \ {\rm\ for\ } n_1<0 {\rm\ or\ } n_2<0,
\eeq
we can calculate many values of $F(m;\, n_1,n_2)$. For example, the sequence $F(2n;\, 0,0)_{n = 0}^\infty$ of
the numbers of lattice walks returning to $(0,0)$ after $2n$ steps starts out as
\begin{eqnarray*}
&& 1,\ 2,\ 11,\ 85,\ 782,\ 8004,\ 88044,\ 1020162,\ 12294260,\ 152787976,\ 1946310467,\ 25302036071, \\
&& 334560525538,\ 4488007049900,\ 60955295750460,\ 836838395382645,\ 11597595644244186,\ \ldots
\end{eqnarray*}
Based on such empirical evidence, Ira Gessel conjectures (cf.\ \cite[p.\ 4]{zeil}) that
\beq
\label{gessel}
F(2n;\,0,0) = 16^n \frac{(1/2)_n(5/6)_n}{(2)_n(5/3)_n}.
\eeq
Therefore we will call the numbers $F(2n;\, 0,0)$ {\em Gessel numbers}.

\section{Similar conjectures}

Analogous conjectures can be made about other points. For example, we conjecture that
\[
F(2n;\,0,1)\ =\ 16^n \frac{\left(1/2\right)_n}{(3)_n}
\left(\frac{5}{27} \frac{\left(7/6\right)_n}
   {\left(7/3\right)_n} + \frac{\left(111 n^2+183 n-50\right)}{270} \frac{
   \left(5/6\right)_n}{\left(8/3\right)_n }\right).
\]
More generally, looking at the points $(0,k)$ it seems that
\[
F(2n;\,0,k)\ =\
16^n \frac{\left(1/2\right)_n }{\left(k+2\right)_n}
\left(\frac{
   \left(7/6\right)_n
   }{\left((3k+4)/3\right)_n
    }p_k(n)
   + \frac{
   \left(5/6\right)_n
   }{\left((3k+5)/3\right)_n
    }q_k(n)\right)
\]
where $p_k(n)$ is a polynomial of degree $2k-2$ and $q_k(n)$ is a polynomial of degree $2k$.
These two polynomial sequences seem to be non-holonomic.

At $(2n+2k;\,0,n)$ we seem to have
\[
\begin{array}{lll}
F(2n;\,0,n) &=& 4^n \frac{\left(3/2\right)_n }{(2n+1) \left(2\right)_n}\ =\ \frac{4^n \left(1/2\right)_n}{(2)_n}, \\
F(2n+2;\,0,n) &=& \frac{2^{2 n+1} (n+1) \left(3/2\right)_n}{(3)_n}, \\
F(2n+4;\,0,n) &=& \frac{4^n (n+1) \left(8 n^2+32 n+33\right) \left(3/2\right)_n}{3\, (4)_n}, \\
F(2n+6;\,0,n) &=& \frac{4^{n-1} (n+1) \left(64 n^4+672 n^3+2648 n^2+4641 n+3060\right)
   \left(3/2\right)_n}{9\, (5)_n}, {\rm\ etc.,}
\end{array}
\]
from which we conjecture that
\beq
\label{vertexcess}
F(2n+2k;\,0,n)\ =\
4^n \frac{\left(3/2\right)_n }{\left(k+2\right)_n}\,r_k(n)
\eeq
where $r_k(n)$ is a polynomial of degree $2k-1$ divisible by $n+1$
for $k \ge 1$, and $r_0(n) = 1/(2n+1)$. It seems that this polynomial sequence is non-holonomic. 

Another empirical observation is that $g(n) = F(2n+1;\,1,0)$ likely satisfies the second-order recurrence
\beq
\label{2ndOrd}
\begin{array}{rlll}
       (n+3) (3n+7) (3n+8)       & \!\!\!g(n+1) && \\
-\ 8\, (2n+3) (18n^2 + 54n + 35) & \!\!\!g(n)   && \\
+\ 256\, n\, (3n+1) (3n+2)       & \!\!\!g(n-1) &=& 0.
\end{array}
\eeq
According to algorithm {\em Hyper}, this recurrence
has no hypergeometric solutions. But $F(2n;\,2,0)$ seems to be non-holonomic.

On the other hand, we seem to have
\[
\begin{array}{lll}
F(n;\,n,0) &=& 1, \\
F(n+2;\,n,0) &=& \frac{1}{2} (n+1) (n+4), \\
F(n+4;\,n,0) &=& \frac{1}{12} (n+1) \left(n^3+15 n^2+74 n+132\right), \\
F(n+6;\,n,0) &=& \frac{1}{144} (n+1) \left(n^5+32 n^4+407 n^3+2620 n^2+8604 n+12240\right), {\rm\ etc.,}
\end{array}
\]
from which we conjecture that
\beq
\label{horexcess}
F(n+2k;\,n,0)\ =\ s_k(n)
\eeq
where $s_k(n)$ is a polynomial of degree $2k$ with leading coefficient $\frac{1}{k!(k+1)!}$, which is 
divisible by $n+1$ when $k \ge 1$. This polynomial sequence seems to be non-holonomic. While $s_k(0) =
F(2k;\,0,0)$ is hypergeometric by Gessel's conjecture, and $s_k(1)$ seems to be holonomic of order two
as per (\ref{2ndOrd}), the sequences $s_k(2)$, $s_k(3)$, $\ldots$ all seem
to be non-holonomic. On the other hand, the coefficient sequences $[n^{2k}]s_k(n)$, 
$[n^{2k-1}]s_k(n)$, $[n^{2k-2}]s_k(n)$, $\ldots$ all seem to be hypergeometric,
again harboring a polynomial sequence of increasing degrees.

\section{Some values of the numbers $F(m;\, n_1,n_2)$}

\begin{Proposition} 
\label{p1}
$F(m;\, n_1,n_2) \ne 0$ only if
\begin{itemize}
\item[{\rm (i)}] $m \equiv n_1 \pmod 2$,
\item[{\rm (ii)}] $n_1 \le m$, 
\item[{\rm (iii)}] $n_2 \le \frac{1}{2}(n_1 + m)$.
\end{itemize}
\end{Proposition}

{\em Proof:} \ As $F(m;\, n_1,n_2) \ne 0$, there exists a walk $w$ from $(0,0)$ to $(n_1,n_2)$ having $m$ steps.
Assume that out of these $m$ steps,
$a$, $b$, $c$ resp.\ $d$ are $E$, $W$, $NE$ resp.\ $SW$ steps. Then
\begin{eqnarray*}
a+b+c+d &=& m,\\
a-b+c-d &=& n_1, \\
c-d &=& n_2,
\end{eqnarray*}
hence $2a+2c = m+n_1$, so $m \equiv n_1 \pmod 2$. Also, $n_1 = a-b+c-d \le a+b+c+d = m$, and
$n_2 = c-d \le a+c = (m+n_1)/2$. \qed

\begin{Theorem}
\label{Catalan}
Let $s(n_1, n_2)$ be the length of a shortest walk $w$ from $(0,0)$ to $(n_1, n_2)$.
\begin{itemize}
\item[{\rm (i)}] If $n_1 \ge n_2$ then $w$ uses $E$, $NE$ steps only, and
\[
F(n_1;\,n_1,n_2)\ =\ s(n_1, n_2)\ =\ {n_1 \choose n_2}.
\]
\item[{\rm (ii)}] If $n_1 \le n_2$ then $w$ uses $W$, $NE$ steps only, and
\[
F(2n_2 - n_1;\,n_1,n_2)\ =\ s(n_1, n_2)\ =\ \frac{n_1+1}{2n_2-n_1+1}{2n_2-n_1+1 \choose n_2 + 1}.
\]
\end{itemize}
\end{Theorem}

{\em Proof:} \ Denote by $a$, $b$, $c$ resp.\ $d$ the numbers of $E$, $W$, $NE$ resp.\ $SW$ steps
in $w$, and by $m$ the length of $w$.
\begin{itemize}
\item[{\rm (i)}] $n_1 \ge n_2$: By Proposition \ref{p1}(ii), $m \ge n_1$. The walk consisting of $n_2$ $NE$ steps
followed by $n_1-n_2$ $E$ steps ends at $(n_1, n_2)$ and has length $n_1$, so it is the shortest such walk. Hence
$s(n_1, n_2) = F(n_1;\,n_1,n_2)$. From $a+b+c+d=m=n_1=a-b+c-d$ it follows that $b=d=0$, so $w$ uses $E$ and $NE$ steps only.
From $a+c=m=n_1$ and $c=n_2$ it follows that $a=n_1-n_2$. Thus $w$ must contain $n_1-n_2$ $E$ steps and $n_2$ $NE$ steps.
There is no restriction on the order of these steps, so there are ${n_1 \choose n_2}$ such walks.

\item[{\rm (ii)}] $n_1 \le n_2$: By Proposition \ref{p1}(iii), $m \ge 2n_2-n_1$. The walk consisting of $n_2$ $NE$ steps
followed by $n_2-n_1$ $W$ steps ends at $(n_1, n_2)$ and has length $2n_2-n_1$, so it is the shortest such walk.
Hence $s(n_1, n_2) = F(2n_2-n_1;\,n_1,n_2)$.
In the inequality $m=a+b+c+d \ge -a+b+c-d = 2(c-d)-(a-b+c-d)=2n_2-n_1$ equality holds when $a=d=0$, so $w$ uses 
$W$ and $NE$ steps only. Conversely, any walk $w$ using $W$ and $NE$ steps only satisfies $-b+c=n_1$ and $c=n_2$,
so $b=n_2-n_1$ and $m=b+c=2n_2-n_1$, implying that $w$ is a shortest such walk. Thus to compute the number 
$s(n_1, n_2)$ of these walks it suffices to enumerate lattice walks from $(0,0)$ to $(n_1, n_2)$ staying
in the octant $0\le n_1\le n_2$ and using $W$ and $NE$ steps only. We have the recurrence
\[
s(n_1, n_2)\ =\ s(n_1+1, n_2) + s(n_1-1, n_2-1) \quad {\rm for\ } n_2 \ge n_1 + 1 \ge 1
\]
with boundary conditions 
\begin{eqnarray*}
s(-1, n_2)  &=& 0, \quad {\rm for\ } n_2 \ge 0, \\
s(n_1, n_1) &=& 1, \quad {\rm for\ } n_1 \ge 0.
\end{eqnarray*}
To this end, we transform the octant $0\le n_1\le n_2$ into the first quadrant by the linear map 
\[
\mu:\ (1,1)\mapsto (1,0), \ (0,1)\mapsto (0,1).
\]
The matrix corresponding to
\[
\mu^{-1}:\ (1,0)\mapsto (1,1), \ (0,1)\mapsto (0,1)
\]
in the standard basis of $\rs^2$ is 
\[
M^{-1}\ =\ \left[
\begin{array}{rr}
1 & 0 \\
1 & 1
\end{array}\right],
\]
so the matrix corresponding to $\mu$ is
\[
M\ =\ \left[
\begin{array}{rr}
1 & 0 \\
-1 & 1
\end{array}\right].
\]
Writing ${\bf n} = (n_1, n_2)$, define
\[
u(n_1, n_2)\ =\ s(M^{-1}{\bf n}) = s(n_1, n_1+n_2).
\]
Then $s(n_1, n_2) = u(M{\bf n}) = u(n_1,n_2-n_1)$, and $u$ satisfies the recurrence
\begin{eqnarray}
u(n_1, n_2) &=& s(n_1, n_1+n_2) \nonumber \\ 
&=& s(n_1+1, n_1+n_2) + s(n_1-1, n_1+n_2-1) \nonumber \\ 
&=& u(n_1+1,n_2-1) + u(n_1-1, n_2), \label{recCat}
\end{eqnarray}
for $n_1 \ge 0 \land n_2 \ge 1$, with boundary conditions 
\begin{eqnarray}
u(-1, n_2)  &=& s(-1, n_2-1)\ =\ 0, \quad {\rm for\ } n_2 \ge 1, \\ \label{bd1}
u(n_1, 0) &=& s(n_1, n_1)\ =\ 1, \quad {\rm for\ } n_1 \ge 0. \label{bd2}
\end{eqnarray}
Let
\[
f(x,y)\ =\ \sum_{n_1,n_2\ge 0} u(n_1, n_2) x^{n_1}y^{n_2}
\]
be the generating function of $u$. From (\ref{recCat}) -- (\ref{bd2}) we obtain
in the usual way the functional equation
\beq
\label{funeq}
(x-x^2-y) f(x,y)\ =\ x - y f(0,y)
\eeq
which can be solved by the kernel method. Since
\[
x^2-x+y\ =\ \left(x-\frac{1+\sqrt{1-4y}}{2}\right)\left(x-\frac{1-\sqrt{1-4y}}{2}\right),
\]
substituting $x = x(y) = \frac{1-\sqrt{1-4y}}{2}$ in (\ref{funeq}) yields
\[
f(0,y)\ =\ \frac{x(y)}{y}\ =\ \frac{1-\sqrt{1-4y}}{2y}\ =\ C(y),
\]
the generating function of Catalan numbers. Hence
\[
f(x,y)\ =\ \frac{x - y\, C(y)}{x-x^2-y}\ =\ -\frac{1}{x-\frac{1+\sqrt{1-4y}}{2}}
\ =\ \frac{C(y)}{1-x\, C(y)}.
\]
Following \cite[p.\ 154]{riordan}, this can be expanded into
\[
f(x,y)\ =\ \sum_{n_1=0}^\infty x^{n_1} C(y)^{n_1+1}
\ =\ \sum_{n_1,n_2\ge 0} \frac{n_1+1}{2n_2+n_1+1}{2n_2+n_1+1 \choose n_2} x^{n_1}y^{n_2},
\]
so we read off
\[
u(n_1, n_2)\ =\ \frac{n_1+1}{2n_2+n_1+1}{2n_2+n_1+1 \choose n_2}
\]
and, finally,
\[
F(2n_2-n_1;\,n_1,n_2)\ =\ s(n_1, n_2)\ =\ u(n_1, n_2-n_1)\ =\ \frac{n_1+1}{2n_2-n_1+1}{2n_2-n_1+1 \choose n_2-n_1}.
\]
\qed
\end{itemize}

\begin{Corollary}
For all $n \ge 0$, we have
\begin{itemize}
\item[{\rm (i)}] $F(n;\,n,0) = 1$,
\item[{\rm (ii)}] $F(2n;\,0,n) = C_n$,
the $n$-th Catalan number.
\end{itemize}
\end{Corollary}

{\em Proof:} \ By Theorem \ref{Catalan}(i), $F(n;\,n,0) = {n \choose 0} = 1$.
By Theorem \ref{Catalan}(ii), 
\[
F(2n;\,0,n) = \frac{1}{2n+1}{2n+1 \choose n+1} = \frac{1}{n+1}{2n \choose n} = C_n.
\] \qed

\section{A functional equation for the generating function}

Let
\beq
\label{G}
G(x,y,z) \ =\ \sum_{m,n_1,n_2 \ge 0}F(m;\,n_1,n_2)\,x^m y^{n_1}z^{n_2}
\eeq
be the generating function of the numbers $F(m;\,n_1,n_2)$. From (\ref{rec}) -- (\ref{bdry})
we obtain in the usual way the following functional equation satisfied by $G(x,y,z)$:
\beq
\label{kernel}
K(x,y,z)G(x,y,z) \ =\ x(1+z)\,G(x,0,z) + x\, G(x,y,0) - x\, G(x,0,0) - y\, z.
\eeq
Here the polynomial
\beq
\label{K}
K(x,y,z) \ =\ x(1+z)(1 + y^2 z) - y\, z
\eeq
is called the {\em kernel\/} of equation (\ref{kernel}).

In order to simplify (\ref{kernel}), we introduce another generating function
\[
H(x,y,z)\ =\ K(x,y,z)G(x,y,z) + y\, z.
\]
Since 
\begin{eqnarray}
H(x,0,z) &=& x(1+z)\, G(x,0,z), \label{Hx0z} \\ 
H(x,y,0) &=& x\, G(x,y,0), \label{Hxy0} \\
H(x,0,0) &=& x\, G(x,0,0), \nonumber
\end{eqnarray}
equation (\ref{kernel}) becomes
\beq
\label{Hkernel}
H(x,y,z)\ =\ H(x,0,z) + H(x,y,0) - H(x,0,0).
\eeq
This is the functional equation that we will work with in the sequel. Write
\[
H(x,y,z)\ =\ \sum_{m,n_1,n_2 \ge 0}\widetilde{F}(m;\,n_1,n_2)\,x^m y^{n_1}z^{n_2}.
\]
It follows from (\ref{Hkernel}), (\ref{Hxy0}) and (\ref{Hx0z}) that
\begin{eqnarray*}
\widetilde{F}(m;\,n_1,n_2) &=& 0 \qquad{\rm if\ } n_1 n_2 \ne 0, \\
\widetilde{F}(m;\,n_1,0)\ &=& F(m - 1;\, n_1, 0), \\
\widetilde{F}(m;\,0,n_2)\ &=& F(m - 1;\, 0, n_2) + F(m - 1;\, 0, n_2 - 1).
\end{eqnarray*}
Note that $F(2n;\, 0, 0) = \widetilde{F}(2n+1;\,0,0)$, and Gessel's conjecture (\ref{gessel}) can be stated as
\[
G(x,0,0)\ =\ \frac{H(x,0,0)}{x} \ =\ 
_3F_2\left(\left.\begin{array}{c}5/6,\, 1/2,\, 1 \\ 2,\, 5/3 \end{array}\right| 16 x^2 \right)
\ =\ \frac{_2F_1\left(\left.\begin{array}{c}-1/2,\, -1/6 \\ 2/3 \end{array}\right| 16 x^2 \right) - 1}{2x^2}.
\]
In analogy to (\ref{vertexcess}) we conjecture that
\[
\widetilde{F}(2n+2k+1;\,0,n)\ =\ 4^n \frac{\left(1/2\right)_n}{\left(k+2\right)_n} \,\widetilde{r}_k(n)
\]
where $\widetilde{r}_k(n)$ is a polynomial of degree $2k+1$, and the sequence of polynomials $\widetilde{r}_k(n)$
is not holonomic.

\section{The kernel method}

Equations (\ref{kernel}) resp.\ (\ref{Hkernel}) cannot be solved right away because they seem to contain 
other unknown functions beside the full generating functions $G(x,y,z)$ resp.\ $H(x,y,z)$ (the additional unknowns
being just sections of $G(x,y,z)$ resp.\ $H(x,y,z)$, of course). To obtain more information, we look for 
roots of the kernel w.r.t.\ one of the variables which are power series in the remaining variables. Substituting
such roots into (\ref{kernel}) resp.\ (\ref{Hkernel}) yields additional equations which are free of the term 
containing the full generating function $G(x,y,z)$ resp.\ $H(x,y,z)$.

In our case, the kernel (\ref{K}) is linear in $x$, and quadratic in $y$ and $z$.
The roots of $K(x,y,z) = 0$ w.r.t.\ $y$ resp.\ $z$ are not power series in $x, z$ resp.\ $x, y$.
But solving $K(x,y,z) = 0$ for $x$ yields
\begin{eqnarray}
x &=& x(y,z) \ =\ \frac{yz}{(1+z)(1+y^2 z)}\ =\ \frac{y}{1-y}\left(\frac{z}{1+z}-\frac{y z}{1+y z}\right)\label{xroot}\\
&=& \sum_{n=1}^\infty \sum_{k=1}^{n} (-1)^{n+1} y^k z^n \nonumber
\end{eqnarray}
which is a power series in $y, z$ satisfying $x(0,0)=0$. So we can substitute it into (\ref{Hkernel}) to obtain
\beq
\label{X}
H(x(y,z),0,z) + H(x(y,z),y,0) - H(x(y,z),0,0) \ =\ y\,z
\eeq
where the rational function $x(y,z)$ is given in (\ref{xroot}). This does not help us find $H(x,y,z)$ (or a
non-trivial section of it) directly. However, (\ref{X}) does determine all the 
coefficients of $H(x,y,0)$ and $H(x,0,z)$, and hence also of $H(x,y,z)$ and $G(x,y,z)$. 
How could we exploit this?

\section{Gessel numbers as determinants}

Expand the left-hand side of (\ref{X}) into power series in $y$ and $z$, and equate the coefficient of
$y^u z^v$ to~0 (except for the coefficient of $y\, z$ which is equated to 1). This yields
the following infinite system of linear equations
for the values of $\widetilde{F}(m;\, n_1,n_2)$ on the planes $n_1=0$ and $n_2=0$:
\begin{eqnarray}
\lefteqn{
\sum_{m, n_2 \ge 0 \atop m \equiv u \!\!\!\!\!\pmod 2}
{-m \choose \frac{u-m}{2}}
{-m \choose v-n_2-\frac{u+m}{2}}
\widetilde{F}(m;\,0,n_2)}  \nonumber \\
&+&
\sum_{m \ge 0, n_1 \ge 1 \atop m+n_1 \equiv u \!\!\!\!\!\pmod 2}
{-m \choose \frac{u-m-n_1}{2}}
{-m \choose v-\frac{u+m-n_1}{2}}
\widetilde{F}(m;\,n_1,0) \label{sys} \\
&=&
\left\{
\begin{array}{ll}
1, & u = v = 1, \\
0, & {\rm otherwise},
\end{array}
\right. \qquad {\rm for\ all\ } u, v \ge 0. \nonumber
\end{eqnarray}
Note that both sums are finite since a binomial coefficient vanishes when its lower symbol is
negative. So $m \le u$ and $n_2 \le v - u/2$ in the first sum, and $m + n_1 \le u$ in the second.
To put this system into a more compact form, pack its unknowns into an infinite matrix 
$[f(i,j)]_{i,j = 0}^\infty$ defined by
\beq\label{fij}
f(i,j) \ =\ \left\{
\begin{array}{ll}
\widetilde{F}(i;\, 0,j-i), & i \le j, \\
\widetilde{F}(j;\, i-j,0), & i \ge j,
\end{array}
\right.
\eeq
or graphically,
\begin{eqnarray*}
f &=&
\left[
\begin{array}{ccccccc}
\widetilde{F}(0;\, 0, 0)& \widetilde{F}(0;\, 0, 1)& \widetilde{F}(0;\, 0, 2) & \widetilde{F}(0;\, 0, 3) & \widetilde{F}(0;\, 0, 4) & \widetilde{F}(0;\, 0, 5) & \ \;\cdots \\
\widetilde{F}(0;\, 1, 0)& \widetilde{F}(1;\, 0, 0)& \widetilde{F}(1;\, 0, 1) & \widetilde{F}(1;\, 0, 2) & \widetilde{F}(1;\, 0, 3) & \widetilde{F}(1;\, 0, 4) & \ \;\cdots \\
\widetilde{F}(0;\, 2, 0)& \widetilde{F}(1;\, 1, 0)& \widetilde{F}(2;\, 0, 0) & \widetilde{F}(2;\, 0, 1) & \widetilde{F}(2;\, 0, 2) & \widetilde{F}(2;\, 0, 3) & \ \;\cdots \\
\widetilde{F}(0;\, 3, 0)& \widetilde{F}(1;\, 2, 0)& \widetilde{F}(2;\, 1, 0) & \widetilde{F}(3;\, 0, 0) & \widetilde{F}(3;\, 0, 1) & \widetilde{F}(3;\, 0, 2) & \ \;\cdots \\
\widetilde{F}(0;\, 4, 0)& \widetilde{F}(1;\, 3, 0)& \widetilde{F}(2;\, 2, 0) & \widetilde{F}(3;\, 1, 0) & \widetilde{F}(4;\, 0, 0) & \widetilde{F}(4;\, 0, 1) & \ \;\cdots \\
\widetilde{F}(0;\, 5, 0)& \widetilde{F}(1;\, 4, 0)& \widetilde{F}(2;\, 3, 0) & \widetilde{F}(3;\, 2, 0) & \widetilde{F}(4;\, 1, 0) & \widetilde{F}(5;\, 0, 0) & \ \;\cdots \\
\vdots                  & \vdots                  & \vdots                   & \vdots                   & \vdots                   & \vdots                   & \ \;\ddots \\
\end{array}
\right] \\
&=&
\left[
\begin{array}{ccccccccccccccc}
      0 & 0 & 0 & 0 & 0 & 0 & 0 & 0 & 0 & 0 & 0 & 0 & 0 & 0 & \cdots \\
      0 & 1 & 1 & 0 & 0 & 0 & 0 & 0 & 0 & 0 & 0 & 0 & 0 & 0 & \cdots \\
      0 & 0 & 0 & 0 & 0 & 0 & 0 & 0 & 0 & 0 & 0 & 0 & 0 & 0 & \cdots \\
      0 & 0 & 1 & 2 & 3 & 1 & 0 & 0 & 0 & 0 & 0 & 0 & 0 & 0 & \cdots \\
      0 & 0 & 0 & 0 & 0 & 0 & 0 & 0 & 0 & 0 & 0 & 0 & 0 & 0 & \cdots \\
      0 & 0 & 0 & 1 & 5 & 11 & 19 & 10 & 2 & 0 & 0 & 0 & 0 & 0 & \cdots \\
      0 & 0 & 0 & 0 & 0 & 0 & 0 & 0 & 0 & 0 & 0 & 0 & 0 & 0 & \cdots \\
      0 & 0 & 0 & 0 & 1 & 9 & 37 & 85 & 158 & 103 & 35 & 5 & 0 & 0 & \cdots \\
      0 & 0 & 0 & 0 & 0 & 0 & 0 & 0 & 0 & 0 & 0 & 0 & 0 & 0 & \cdots \\
      0 & 0 & 0 & 0 & 0 & 1 & 14 & 87 & 332 & 782 & 1521 & 1126 & 499 & 126 & \cdots \\
      0 & 0 & 0 & 0 & 0 & 0 & 0 & 0 & 0 & 0 & 0 & 0 & 0 & 0 & \cdots \\
      0 & 0 & 0 & 0 & 0 & 0 & 1 & 20 & 172 & 911 & 3343 & 8004 & 16056 & 12941 & \cdots \\
      0 & 0 & 0 & 0 & 0 & 0 & 0 & 0 & 0 & 0 & 0 & 0 & 0 & 0 & \cdots \\
      0 & 0 & 0 & 0 & 0 & 0 & 0 & 1 & 27 & 305 & 2096 & 10147 & 36350 & 88044 & \cdots \\
\vdots & \vdots & \vdots & \vdots & \vdots & \vdots & \vdots & \vdots & \vdots & \vdots & \vdots & \vdots & 
\vdots & \vdots & \ddots \\
\end{array}
\right].
\end{eqnarray*}
Inverting this transformation, we clearly have
\begin{eqnarray*}
\widetilde{F}(m;\,0,n_2) &=& f(m, m+n_2), \\
\widetilde{F}(m;\,n_1,0) &=& f(m+n_1, m).
\end{eqnarray*}
Using this in (\ref{sys}), together with the change of variables $i=m$, $j=m+n_2$ in the first sum, and
$i=m+n_1$, $j=m$ in the second, the left-hand side of (\ref{sys}) changes into
\[
\sum_{j \ge i \ge 0 \atop i \equiv u \!\!\!\!\!\pmod 2}
{-i \choose \frac{u-i}{2}}
{-i \choose v-j-\frac{u-i}{2}}
f(i,j) \ +
\sum_{i \ge j+1 \ge 1 \atop i \equiv u \!\!\!\!\!\pmod 2}
{-j \choose \frac{u-i}{2}}
{-j \choose v-j-\frac{u-i}{2}}
f(i,j) \]
which allows us to combine the two sums into a single one, and so to rewrite (\ref{sys}) as 
\beq
\label{newsys}
\sum_{i, j \ge 0 \atop i \equiv u \!\!\!\!\!\pmod 2}
{-\min \{i,j\} \choose \frac{u-i}{2}}
{-\min \{i,j\} \choose v-j-\frac{u-i}{2}}
\,f(i,j)  
\ =\ 
\left\{
\begin{array}{ll}
1, & u = v = 1, \\
0, & {\rm otherwise},
\end{array}
\right. \ {\rm for\ all\ } u, v \ge 0.
\eeq
Denote the above equation by $E(u,v)$, and let $c(u,v,i,j)$ be the coefficient of $f(i,j)$ in $E(u,v)$: 
\beq
\label{c}
c(u,v,i,j)\ =\
\left\{
\begin{array}{cl}
\displaystyle {-\min \{i,j\} \choose \frac{u-i}{2}}
{-\min \{i,j\} \choose v-j-\frac{u-i}{2}}, & i \equiv u \!\!\!\!\!\pmod 2, \\
0, & {\rm otherwise}. 
\end{array}
\right.
\eeq
\begin{Proposition}
\label{ord}
\ Let $u, v, i, j \ge 0$. Then
\begin{itemize}
\item[\rm (i)] $c(u,v,u,v) = 1$,
\item[\rm (ii)] $c(u,v,i,j) = 0$ if $i > u$ or $j > v$.
\end{itemize}
\end{Proposition}

{\em Proof:\/}
Assertion (i) is clear from (\ref{c}). To prove (ii), assume that $i > u$ or $j > v$.
If $i > u$ then $(u-i)/2 < 0$. Otherwise $i \le u$. Then, by assumption, $j > v$, and so $v-j-(u-i)/2 < 0$.
In either case, $c(u,v,i,j) = 0$. \qed

\

\noindent
Proposition \ref{ord} implies that we can compute $f(u,v)$ from $E(u,v)$, provided that we have already computed
$f(i,j)$ for all $(i,j) \ne (u,v)$ such that $i \le u$ and $j \le v$. In other words, the system (\ref{newsys})
is a linear recurrence from which all the $f(i,j)$ can be computed one by one, in any order compatible with the
standard componentwise partial order on $\ns \times \ns$. 
Nevertheless, we'll continue to regard (\ref{newsys}) as an infinite system of linear equations, and will rewrite it
in the form ${\bm A}{\bm x} = {\bm b}$ where ${\bm A}$ is an infinite matrix and ${\bm x},{\bm b}$ are infinite
vectors. Then we can rephrase Proposition \ref{ord} as follows:

\begin{Corollary}
\label{corord}
Let $\rho: \ns\times\ns \rightarrow \ns$ be a monotonic bijection in the sense that $\rho(a,b) \le \rho(c,d)$
whenever $a \le c\ \land\ b \le d$. Define ${\bm A} = [a(n,k)]_{n,k=0}^\infty$, ${\bm x} = [x(k)]_{k=0}^\infty$
and ${\bm b} = [b(n)]_{n=0}^\infty$ by
\begin{eqnarray}
a(n,k) &=& c(u,v,i,j), \label{ank} \\ 
x(k) &=& f(i,j), \nonumber \\
b(n) &=& \left\{
\begin{array}{ll}
1, & u = v = 1, \\
0, & {\rm otherwise},
\end{array}
\right. \nonumber
\end{eqnarray}
where $c$ resp.\ $f$ is given by (\ref{c}) resp.\ (\ref{fij}), $(u,v) = \rho^{-1}(n)$, and
$(i,j) = \rho^{-1}(k)$. Then
\begin{itemize}
\item[\rm (i)] ${\bm A}$ is lower-triangular with unit diagonal,
\item[\rm (ii)] ${\bm A}{\bm x} = {\bm b}$.
\end{itemize}
\end{Corollary}

{\em Proof:\/}
\begin{itemize}
\item[\rm (i)] By Proposition \ref{ord}(i), $a(n,n) = c(u,v,u,v) = 1$, proving that ${\bm A}$ has unit diagonal. 
Now assume that $n < k$. Then $\rho(u,v) < \rho(i,j)$, hence by monotonicity of $\rho$, $i > u$ or $j > v$.
By Proposition \ref{ord}(ii), $a(n,k) = c(u,v,i,j) = 0$, proving that ${\bm A}$ is lower-triangular.
\item[\rm (ii)] Let $n \in \ns$ be arbitrary, and $(u,v) = \rho^{-1}(n)$. Then by (i), the sum 
$\sum_{k=0}^\infty a(n,k)\, x(k)$ exists, and by (\ref{c}) and (\ref{newsys}), 
\[
\sum_{k=0}^\infty a(n,k)\, x(k) = \sum_{i,j\ge 0} c(u,v,i,j)\,f(i,j) = \left\{
\begin{array}{ll}
1, & u = v = 1, \\
0, & {\rm otherwise},
\end{array}
\right.
 =\,\ b(n),
\]
proving that ${\bm A}{\bm x} = {\bm b}$. \qed
\end{itemize}
We can compute a particular component $x(k)$ of the solution vector ${\bm x}$ from the finite 
lower-triangular system with unit diagonal
\[
{\bm A}^{(k)}{\bm x}^{(k)} \ = \ {\bm b}^{(k)}
\]
where
\begin{eqnarray*}
{\bm A}^{(k)} &=& \left[a(i,j)\right]_{i,j=0}^k, \\
{\bm x}^{(k)} &=& \left(x(j)\right)_{j=0}^k, \\
{\bm b}^{(k)} &=& \left(b(i)\right)_{i=0}^k.
\end{eqnarray*}
By Cramer's rule,
\[
x(k)\ =\ \frac{\det \widetilde{\bm A}^{(k)}}{\det {\bm A}^{(k)}}
\ =\ \det \widetilde{\bm A}^{(k)}
\]
where $\widetilde{\bm A}^{(k)}$ is obtained from ${\bm A}^{(k)}$ by replacing its last column with
${\bm b}^{(k)}$.

If $k < \rho(1,1)$ then ${\bm b}^{(k)} = {\bm 0}$ and $x(k)=0$. If $k \ge \rho(1,1)$ then the last column
of $\widetilde{\bm A}^{(k)}$ has a single non-zero entry, 1, at position $\rho(1,1)$. Developing 
$\det \widetilde{\bm A}^{(k)}$ w.r.t.\ this column yields 
\[
x(k)\ =\ \det \widetilde{\bm A}^{(k)}\ =\ \det {\bm T}^{(k)} \det {\bm H}^{(k)}\ =\ \det {\bm H}^{(k)}
\]
where ${\bm T}^{(k)}$ is a $\rho(1,1) \times \rho(1,1)$ lower-triangular matrix with unit diagonal, and 
${\bm H}^{(k)}$ is the $(k-\rho(1,1)) \times (k-\rho(1,1))$ lower-Hessenberg matrix with unit superdiagonal,
composed of the elements in rows $\rho(1,1)$ to $k$ and columns $\rho(1,1)-1$ to $k-1$ of ${\bm A}^{(k)}$.
Thus for Gessel numbers we have
\[
F(2n;\,0,0)\ =\ \widetilde{F}(2n+1;\,0,0)\ =\ f(2n+1,2n+1)\ =\ x(\rho(2n+1,2n+1))\ =\ \det {\bm H}^{(\rho(2n+1,2n+1))}
\]
when $n \ge 0$.
For example, if we use diagonal ordering to pack the unknown $f(i,j)$ into the vector $\bm x$, then
\[
\rho(i,j) \ =\ {i+j+1 \choose 2} + j,
\]
so $\rho(1,1) = 4$ and $\rho(3,3) = 24$. Hence, if $n=1$, $F(2,0,0) = \det {\bm H}^{(24)} = 2$ where
\[
{\bm H}^{(24)}  \ =\ 
\left[
\begin{array}{cccccccccccccccccccc}
      0 & 
    1 & 0 & 0 & 0
       & 0 & 0 & 0 & 0 & 0 & 0 & 0 & 0 & 0 & 0 & 0 & 0 & 0 & 0 &
         0 \\
      0 & 0 & 1 & 0 & 0 & 0 & 0 & 0 & 0 & 0 & 0 & 0 & 0 & 
        0 & 0 & 0 & 0 & 0 & 0 & 0 \\
      0 & 0 & 0 & 1 & 0
       & 0 & 0 & 0 & 0 & 0 & 0 & 0 & 0 & 0 & 0 & 0 & 0 & 0 & 0 &
         0 \\
      -1 & 0 & 0 & 0 & 1 & 0 & 0 & 0 & 0 & 0 & 0 & 0 & 0
         & 0 & 0 & 0 & 0 & 0 & 0 & 0 \\
      0 & 0 & 0 & 0 &
       0 & 1 & 0 & 0 & 0 & 0 & 0 & 0 & 0 & 0 & 0 & 0 & 0 & 0 & 
        0 & 0 \\
      0 & 0 & 0 & 0 & 0 & 0 & 1 & 0 & 0 & 0 & 0 & 0 & 0 & 0 & 
        0 & 0 & 0 & 0 & 0 & 0 \\
      0 & 0 & 0 & 0 &
       0 & 0 & 0 & 1 & 0 & 0 & 0 & 0 & 0 & 0 & 0 & 0 & 0 & 0
         & 0 & 0 \\
      0 & 0 & 0 & -1 & 0 & 0 & 0 & 0 & 1 & 0 & 0 & 0 & 0 & 
        0 & 0 & 0 & 0 & 0 & 0 & 0 \\
      1 & 0 & 0 & 0 & -1 & 0 & 0 & 0 & 0 & 1 & 0 & 0 & 0 & 
        0 & 0 & 0 & 0 & 0 & 0 & 0 \\
      0 & 0 & 0 & 0 & 0 & 0 & 0 & 0 & 0 & 0 & 1 & 0 & 0 & 0 & 
        0 & 0 & 0 & 0 & 0 & 0 \\
      0 & 0 & 0 & 0 & 0 & 0 & 0 & 0 & 0 & 0 & 0 & 1 & 0 & 0 & 
        0 & 0 & 0 & 0 & 0 & 0 \\
      0 & 0 & 0 & 0 & 0 & 0 & 0 & 0 & 0 & 0 & 0 & 0 & 1 & 0 & 
        0 & 0 & 0 & 0 & 0 & 0 \\
      -1 & 0 & 0 & 0 & 0 & 0 & 0 & -1 & 0 & 0 & 0 & 0 & 0 & 
        1 & 0 & 0 & 0 & 0 & 0 & 0 \\
      0 & 0 & 0 & 1 & 0 & 0 & 0 & 0 & -2 & 0 & 0 & 0 & 
        0 & 0 & 1 & 0 & 0 & 0 & 0 & 0 \\
      -1 & 0 & 0 & 0 & 1 & 0 & 0 & 0 & 0 & -1 & 0 & 0 & 0 & 
      0 & 0 & 1 & 0 & 0 & 0 & 0 \\
      0 & 0 & 0 & 0 & 0 & 0 & 0 & 0 & 0 & 0 & 0 & 0 & 0 & 0
         & 0 & 0 & 1 & 0 & 0 & 0 \\
      0 & 0 & 0 & 0 & 0 & 0 & 0 & 0 & 0 & 0 & 0 & 0 & 0 & 0
         & 0 & 0 & 0 & 1 & 0 & 0 \\
      0 & 0 & 0 & 0 & 0 & 0 & 0 & 0 & 0 & 0 & 0 & 0 & 0 & 0
         & 0 & 0 & 0 & 0 & 1 & 0 \\
      0 & 0 & 0 & -1 & 0 & 0 & 0 & 0 & 0 & 0 & 0 & 
        0 & -1 & 0 & 0 & 0 & 0 & 0 & 0 & 1 \\
      1 & 0 & 0 & 0 & -1 & 0 & 0 & 1 & 0 & 0 & 0 & 0 & 0
       & -2 & 0 & 0 & 0 & 0 & 0 & 0
\end{array}
\right].
\]

\section{Gessel numbers as multiple sums}

For any (finite or infinite) lower-triangular matrix ${\bm A} = [a(k,m)]_{k,m\ge 0}$ with unit diagonal entries,
the lower-triangle elements of its inverse ${\bm A}^{-1} = [\bar{a}(k,m)]_{k,m\ge 0}$ are given by
\beq
\label{bara}
\bar{a}(k,m)\ =\ \sum_{j=1}^{k-m} (-1)^j \sum_{m=\lambda_0 < \lambda_1 < \cdots < \lambda_j = k}
\ \prod_{i=1}^j a(\lambda_i, \lambda_{i-1})
\eeq
when $k > m$. Therefore for any ordering $\rho$ as described in Corollary \ref{corord},
\[
F(2n;\,0,0)\ =\ \bar{a}(\rho(2n+1,2n+1),\rho(1,1))
\]
where $\bar{a}, a$, and $c$ are given in (\ref{bara}), (\ref{ank}), and (\ref{c}),
respectively.

\section{The solution vector}

Here we describe a few properties of the solution vector $\mathbf{x}^{(n)}$.
 
 We'll think of the $n^2\times 1$ vector $\mathbf{x}^{(n)}$ as consisting of the concatenation of $\lfloor{n^2/(2n+1)}\rfloor$ vectors $\mathbf{u}_i$ $(i=1,2,\dots)$, each of length $2n+1$, plus one more, of length $n^2$ mod $2n+1$ . Each of these vectors $\mathbf{u}_i$, except possibly the last, consists first of a certain universal sequence of length $2i$, followed by $(2n+1-2i)$ 0's. The last one consists of as much of the next universal sequence as there is room for. As $n$ increases these vectors remain unchanged, and a new one appears at the end. Each of these universal sequences ends in a Catalan number.

The first several of these universal sequences are
\begin{verbatim}
1, 1
1, 2, 3, 1
1, 5, 11, 19, 10, 2
1, 9, 37, 85, 158, 103, 35, 5
1, 14, 87, 332, 782, 1521, 1126, 499, 126, 14
1, 20, 172, 911, 3343, 8004, 16056, 12941, 6765, 2296, 462, 42
1, 27, 305, 2096, 10147, 36350, 88044, 180621, 154750, 
    90681, 37178, 10254, 1716, 132
1, 35, 501, 4300, 25927, 118472, 417565, 1020162, 
    2128824, 1910006, 1217523, 570409, 193137, 44913, 6435, 429
\end{verbatim}
The sequence of next-to-last members of the above is also holonomic, but the third-from-last sequence might or might not be.

Can we find these sequences by dealing only with the corresponding sections of the matrix?


\begin{thebibliography}{9}
\bibitem{zeil}
M.\ Kauers, D.\ Zeilberger, The quasi-holonomic ansatz and restricted lattice walks, to appear in
{\em J.\ Difference Equations and Applications}.
\bibitem{riordan}
J.\ Riordan, {\em Combinatorial Identities}, John Wiley \& Sons, Inc., New York-London-Sydney 1968.
\end{thebibliography}
\end{document}